\documentclass{amsart}

\usepackage{amsmath}
\usepackage{amsfonts}
\usepackage{amssymb,amscd,amsthm}
\usepackage{geometry}
\usepackage{mathrsfs}
\usepackage[bookmarksnumbered, colorlinks, linktocpage, plainpages]{hyperref}
\usepackage{nameref}
\usepackage{wasysym,amssymb,eufrak,indentfirst,color}
\usepackage{extarrows}

\usepackage{enumitem}

\usepackage{mathtools}
 \usepackage{appendix}
\usepackage{extarrows}
\usepackage{setspace}

\usepackage{hyperref}
\usepackage{cleveref}

\usepackage{cases}
\usepackage{graphicx}
\usepackage{ragged2e}

\newtheorem{theorem}{Theorem}
 \newtheorem{proposition}{Proposition}
 \newtheorem{definition}{Definition}
 \newtheorem{remark}{Remark}
  \newtheorem{example}{Example}

\newtheorem{lemma}{Lemma}

\numberwithin{equation}{section}
\newcommand{\R}{\mathbb{R}}

\def\R{\mathbb{R}}

\numberwithin{equation}{section}

\def\R{\mathbb{R}}

\begin{document}

\title[The Twisted Group Algebra Structure of the Cayley-Dickson Algebra]{The Twisted Group Algebra Structure of the Cayley-Dickson Algebra}

\author{Guangbin Ren}
\email[G.~Ren]{rengb@ustc.edu.cn}
\address{Department of Mathematics, University of Science and Technology of China, Hefei 230026, China}

\author{Xin Zhao}
\email[X.~Zhao]{zx130781@mail.ustc.edu.cn}
\address{Department of Mathematics, University of Science and Technology of China, Hefei 230026, China}

 \date{\today}


\thanks{This work was supported by the NNSF of China (12171448).}

\maketitle
\begin{abstract}
The Cayley-Dickson algebra has long been a challenge due to the lack of an explicit multiplication table. Despite being constructible through inductive construction, its explicit structure has remained elusive until now. In this article, we propose a solution to this long-standing problem by revealing the Cayley-Dickson algebra as a twisted group algebra with an explicit twist function $\sigma(A,B)$. We show that this function satisfies the equation $$e_Ae_B=(-1)^{\sigma(A,B)}e_{A\oplus B}$$ and provide a formula for the relationship between the Cayley-Dickson algebra and split Cayley-Dickson algebra, thereby giving an explicit expression for the twist function of the split Cayley-Dickson algebra. Our approach not only resolves the lack of explicit structure for the Cayley-Dickson algebra and split Cayley-Dickson algebra but also sheds light on the algebraic structure underlying this fundamental mathematical object.
\end{abstract}

    \tableofcontents

\section{Introduction}

The Cayley-Dickson algebra has been a fundamental mathematical structure for many years due to its applications in mathematics and physics \cite{SCh1995, HR2022, MB2019, Sw2010,Gr2021, BK1998,Mr2021,Mr2023,MY2014}. However, a major challenge associated with this algebra is the lack of an explicit multiplication table.

Recent research has approached this problem by treating the Cayley-Dickson algebra as a twisted group algebra \cite{AM1999,BS1970,RZ2022,R1971}. In this framework, the twist function $\sigma(A,B)$ plays a crucial role in determining the multiplication table. Although explicit twist functions exist for certain cases, the general case has remained elusive, with the twist function thought to be highly chaotic.

Previous work has attempted to address this problem through the quasi-algebra approach. For instance, Albuquerque and Majid \cite{AM1999} provided recurrence formulas, and Flaut and Boboescu \cite{FB2021} proposed an alternative method. Nevertheless, an explicit formula for the twist function has yet to be discovered.

In this article, we present an explicit formula for the twist function of the Cayley-Dickson algebra.  This algebra is a generalization of the complex numbers obtained by repeatedly applying the Cayley-Dickson process. The algebra is denoted by $\mathbb D^{(n)}$ and has dimension $2^n$, where $n$ is a non-negative integer. In the multiplication table of $\mathbb D^{(n)}$, the product of two basis elements $e_A$ and $e_B$ is given by:
$$e_Ae_B=(-1)^{\sigma(A,B)}e_{A\oplus B},$$
where $A$ and $B$ are binary sequences of length $n$ that represent the basis elements, and $\oplus$ denotes bitwise addition.

The twist function $\sigma(A,B)$ determines the sign of the product $e_Ae_B$ and depends on the degrees and values of $A$ and $B$. Theorem \ref{0515-1855} provides an explicit formula for $\sigma(A,B)$:
 \begin{eqnarray}\label{0502-2056} \sigma(A,B)=  \left\{
 \begin{aligned} & 0, \quad  & &AB=0,
\\
&1, \quad  & &A=B\neq 0,
\\
       & b_l+\sum\limits_{i=l}^{n-1}\phi(a_i,b_i),\quad  & & A\neq B, \quad AB\neq 0,\quad \deg A> \deg B,\\
       & 1+a_l+\sum\limits_{i=l}^{n-1}\phi(a_i,b_i),\quad & &A\neq B, \quad AB\neq 0, \quad  \deg A<\deg B,\\
       & a_l+\sum\limits_{i=l}^{n-1}\phi(a_i,b_i),\quad  & &A\neq B,  \quad AB\neq 0,\quad \deg A=\deg B,
            \end{aligned}  \right.
     \end{eqnarray}
Here we write $A$ and $B$ in the form of binary sequences. 
 \begin{eqnarray}\label{eq:binary-490-23} A&=&a_02^0+\dots+a_{n-1}2^{n-1},
 \\ \label{eq:binary-491-23}
     B&=&b_02^0+\dots+b_{n-1}2^{n-1}.
     \end{eqnarray} 
Each coefficient $a_j, b_j$ is either 0 or 1.
 We define $\deg\, A$  to be the index of the first non-zero term in the binary sequence representing $A$. Namely, 
$$\deg\, A=\min\{i\geqslant 0: a_i\neq 0 \}.$$
 $\deg B$ can be defined similarly.
The function $l=l(A, B)$ is  defined by 
 \begin{eqnarray}\label{eq:phi-939} l=\max\{\deg A, \, \deg B, \, \deg (A\oplus B)\}.
\end{eqnarray}
and 
 the function $\phi(a,b)$ is defined as 
 \begin{eqnarray}\label{eq:binary-490-99} \phi(a,b)=\begin{cases}
0 & \text{if } a=b=0, \\
1 & \text{otherwise}.
\end{cases}\end{eqnarray}

Using the same method, we can derive similar equations for split Cayley-Dickson algebras. These algebras are hyperbolic counterparts of Cayley-Dickson algebras and have a wide range of applications, which are discussed in detail in \cite{Gr2021}.

\bigskip

\section{Cayley-Dickson algebras}

In this section, we will delve into the internal structure of Cayley-Dickson algebras, which are defined recursively using the Cayley-Dickson process over $\mathbb{R}$, as explained in \cite{B2002, Gu2020,CR2014,CR2018}. Specifically, this process involves adding one generator to $\mathbb{D}^{(n)}$ to construct $\mathbb{D}^{(n+1)}$.

To define the Cayley-Dickson algebras over $\mathbb{R}$, we start with $\mathbb{D}^{(0)}=\mathbb{R}$ and $\mathbb{D}^{(n)}$, which is an algebra with standard algebraic generators $g_0, g_1, \ldots, g_{n-1}$ satisfying $$g_i^2=-1, \qquad \overline{g_i}=-g_i$$ for $i=0,1,\ldots,n-1$. The algebra $\mathbb{D}^{(n+1)}$ is then obtained from $\mathbb{D}^{(n)}$ by adding a new imaginary unit $g_n$ and applying the multiplication rule 
\begin{equation}\label{eq:CD399} (a+b g_{n})(c+ d g_{n})=(ac-\overline d b)+ (da+b\overline c) g_{n}  \end{equation}
   for any $a,b,c,d\in \mathbb{D}^{(n)}$. This procedure is known as a doubling procedure of a smashed product, as introduced in \cite{SCh1995,S1954}.

The multiplication relations in \eqref{eq:CD399} can be expressed in the following way:   For any $a,b\in \mathbb{D}^{(n)}$, we have
 \begin{eqnarray}\label{CDrelation1}
    g_{n}^2&=&-1,
\\
\label{CDrelation2}
    \bar g_{n}&=&-g_{n}
\\
\label{CDrelation3}
    a(g_{n}b)&=&g_{n}(\bar ab),
\\
\label{CDrelation4}
    (ag_{n})b&=&(a\bar{b})g_{n},
\\
\label{CDrelation5}
    (g_{n}a)(bg_{n})&=&-\overline{ab}.
\end{eqnarray}

There are certain relations that arise from this construction.
\begin{lemma}
    The identities \eqref{CDrelation1}$\sim$\eqref{CDrelation5}  imply the following relations:
    \begin{eqnarray}\label{CDrelationship11}
        ag_{n}&=&g_{n}\bar a,
 \\ \label{CDrelationship12}
        g_{n}a &=& \bar a g_{n},
  \\ \label{CDrelationship13}
        (g_{n}a)b&=&g_{n}(ba),
   \\ \label{CDrelationship14}
        \qquad a(bg_{n})&=&(ba)g_{n},
  \\
  \label{CDrelationship15}
        (g_{n}a)(g_{n}b)&=&-b\bar a,
  \\
  \label{CDrelationship16}
        (ag_{n})(bg_{n})&=&-\bar ba.
    \end{eqnarray}
\end{lemma}
%\begin{remark}
%    $\mathbb{D}^{(1)}=\CC$, complex number,associative and commutative;\\
%    $\mathbb{D}^{(2)}=\mathbb{H}$, quaternion number, associative but noncommutative;\\
%    $\mathbb{D}^{(3)}=\mathbb{O}$, octonions number, nonassociative and noncommutative.\\
%    \D is nonassociative and noncommutative for any $n\geq 3$.
%\end{remark}

% \begin{proposition}\label{prop:basic}
%     In the nonassociative algebra \D, the product of any finite generators only differs by -1 at most.
% \end{proposition}

{\bf Convention:} Multiplication of octonions is conventionally performed from left to right. That is,
\begin{equation}
g_{i_0}g_{i_1}\dots g_{i_k}:=\Big(\big((g_{i_0}g_{i_1})g_{i_2}\big)\dots g_{i_k}\Big).
\end{equation}

\begin{definition} Let $A$
  be a nonnegative integer  with binary expansion $$A=a_02^0+a_12^1+\dots+a_{n-1}2^{n-1},$$ where all $a_i=0,1$. We define
    \begin{equation}\label{def:eA}
       e_A:=g_0^{a_0}g_1^{a_1}\dots g_{n-1}^{a_{n-1}}.
    \end{equation}

\end{definition}

\bigskip

\section{Twisted group algebras}

In this section, we will treat the Cayley-Dickson algebras as group algebras (see \cite{AM1999}) and recall the recursive formula for the corresponding twist functions.

\begin{definition}
    A twisted group algebra $k_FG$ for a finite group G over a field k is a  k-algebra which has a basis$ \{a_g:g\in G\}$ with
    $$a_g\cdot a_{g^\prime} :=F(a_g,a_{g^\prime}) a_{gg^\prime},$$ where
    $F(a_g,a_{g^\prime})\in k\backslash 0.$
\end{definition}

\begin{lemma}   \cite[Prop 4.1]{AM1999} \label{0503-933}
    There exists a 2-cochain map $$F:\mathbb Z_2^n\times \mathbb Z_2^n\longrightarrow \mathbb Z_2,$$
    such that
     $$\mathbb{D}^{(n)}\cong k_FG,$$
    where $$k=\mathbb Z_2, \qquad G=\mathbb Z_2^n.$$
\end{lemma}

\begin{definition}
  Let  $A,B\in \mathbb N$ be  with binary decompositions
    \begin{equation}\label{KAKB}
        \begin{aligned}
          A&=a_02^0+\dots+a_{n-1}2^{n-1},\\
          B&=b_02^0+\dots+b_{n-1}2^{n-1}
        \end{aligned}
    \end{equation}
  with all $a_i,b_i\in \mathbb Z_2$. We define the binary addition as $$A\oplus B:=(a_0+ b_0)2^n+\dots+(a_{n-1}+b_{n-1})2^{n-1}.$$
  Here we make a  convention that all the coefficients
  $$a_i+b_i\in\mathbb Z_2.$$

\end{definition}

 \begin{lemma} For any $A$ and $B$ in $\mathbb N$, there exists a map
 $$\sigma: \mathbb Z\times \mathbb Z \longrightarrow \mathbb Z_2 $$
 such that
\begin{equation}\label{eq:main-formula-293}
e_Ae_B=(-1)^{\sigma(A,B)} e_{A\oplus B}.
\end{equation}
\end{lemma}

\begin{proof}  Denote $$F(A,B)=(-1)^{\sigma(A,B)}.$$
The result follows directly from   \Cref{0503-933}.
\end{proof}

 We will provide the proof of the recursion formula for the twist function $\sigma(A,B)$, as given in \cite{AM1999}, to ensure completeness.

\begin{proposition}\label{0504-1816}
    For any integers  $0\leqslant A,B\leqslant 2^n-1$ and $a_n, b_n\in\{0, 1\}$, we have the  recursion formula
    \begin{equation}\label{1227-2}
        \begin{aligned}
            \sigma(A+a_{n}2^{n},B+b_{n}2^{n})&=\sigma(A,B)(1+b_{n})+\sigma(B,A)b_{n}+\sigma(B,B)a_{n}+a_{n}b_{n}.
        \end{aligned}
     \end{equation}

\end{proposition}
\begin{proof}
  Starting with \eqref{CDrelation4}, we have $e_{A+2^n}e_B=(e_A\bar e_B)e_{2^n}$, which gives us
\begin{equation}\label{0502-1215}
\sigma(A+2^n,B)=\sigma(A,B)+\sigma(B,B).
\end{equation}
Using \eqref{CDrelationship14}, we get
$$e_Ae_{B+2^n}=(e_Be_A)e_{2^n},$$
which leads to
\begin{equation}\label{0502-1216}
\sigma(A,B+2^n)=\sigma(B,A).
\end{equation}
Using \eqref{CDrelationship16}, we have
$$e_{A+2^n}e_{B+2^n}=-\bar e_Be_A,$$
which implies
\begin{equation}\label{0502-1220}
\sigma(A+2^n,B+2^n)=1+\sigma(B,B)+\sigma(B,A).
\end{equation}
Combining \eqref{0502-1215}, \eqref{0502-1216}, and \eqref{0502-1220}, we obtain
$$\sigma(A+a_n2^n,B+b_n2^n)=\sigma(A,B)(1+b_{n})+\sigma(B,A)b_{n}+a_{n}\sigma(B,B)+a_{n}b_{n}.$$
This completes the proof.
\end{proof}

\section{The main result}

This section will present a concrete formula for the twist function of the Cayley-Dickson algebra based on its multiplication table.

\begin{definition} The degree of
   any positive integer $A$
      with the binomial expression
      $$ A=a_02^0+a_12^1+\dots+a_n2^n$$
      is defined to be the index of of the leading term,  i.e.,
    $$\deg\, A:=\min\{i\geqslant 0: a_i\neq 0   \}.$$
\end{definition}

\begin{definition} We define an auxiliary function    $$\phi :\mathbb Z_2\times \mathbb Z_2\longrightarrow \mathbb Z_2,$$ defined by
    $$\phi(a,b)=a+b+ab.$$
\end{definition}

For any given  positive integera $A$ and $B$, 
we introduce a constant
  $$l=l(A, B):=\max\{\deg A,\deg B, \deg (A\oplus B)\}.$$

To establish our main result, we require a crucial lemma.

\begin{lemma}\label{0515-1832}  Let $0\leqslant A, B\leqslant 2^n-1$ be positive integers with binary decompositions
   $$ A=a_02^0+\dots+a_{n-1}2^{n-1},$$
    $$B=b_02^0+\dots+b_{n-1}2^{n-1}.$$

(1)  If $\deg A= \deg B,$ then $$a_l+b_l=1.$$

(2)  We have $$\phi(a_l,b_l)=1.$$
\end{lemma}

\begin{proof}
(1) We denote $$n_0:=\deg A=\deg B.$$
Let $N$  be the  integer  such that
 $$N>n_0$$
 and
$$ \left \{
\begin{array}{lcl} a_{n_0} &=& b_{n_0}
\\
 a_{n_0+1}&=&b_{n_0+1}
 \\
  \qquad &\ldots&
  \\ a_{N-1}&=&b_{N-1}
  \\
    a_{N}&\neq& b_{N}.
\end{array}
\right.
$$
In the finite field $\mathbb Z_2$, this can be rewritten as
 $$ \left \{
\begin{array}{lclcl} a_{n_0}&+& b_{n_0}  &=& 0
\\
 a_{n_0+1}&+&b_{n_0+1} &=& 0
 \\
  \qquad && &\ldots&
  \\ a_{N-1}&+&b_{N-1} &=& 0
  \\
    a_{N}&+& b_{N} &=& 1.
\end{array}
\right.
$$

Consequently, we can rewrite  $A$ and $B$  as
\begin{eqnarray*}
A &=&a_{n_0}2^{n_0}+a_{n_0+1} 2^{n_0+1}+\cdots+ a_{N-1} 2^{N-1}+ a_{N} 2^{N}+\cdots,
\\
B &=&a_{n_0}2^{n_0}+a_{n_0+1} 2^{n_0+1}+\cdots+ a_{N-1} 2^{N-1}+ b_{N} 2^{N}+\cdots.
\end{eqnarray*}
Since $a_{N} \neq  b_{N}$, it follows from the definition of the binary addition that
\begin{eqnarray*} A\oplus B&=&(a_N+b_N)2^N+\cdots
\\
&=&2^N+\cdots. \end{eqnarray*}
  This means
$$\deg(A\oplus B)=N$$
 Consequently,  $$a_l+b_l=1.$$

(2) It is clear that  $$\phi(a,b)=0$$
if and only if$$ a=b=0.$$
Since
$$a_l+b_l=1$$
in the case where $\deg A=\deg B$ due to the just proved result in  (1), it follows that $$\phi(a_l,b_l)=1.$$
For the remaining case where  $\deg A\neq \deg B$,
we apply the definition of $l$ to deduce that
 $$a_l=1, \qquad \mbox{if}\,\, \deg A>\deg B,$$
$$b_l=1, \qquad  \mbox{if}\,\, \deg B>\deg A.$$
In both cases, we always have  $$\phi(a_l,b_l)=1.$$
\end{proof}

 Now we arrive at our main result.

\begin{theorem}\label{0515-1855} In the multiplicatin table of the Cayley-Dickson algebra $\mathbb D^{(n)}$ expressed as $$e_Ae_B=(-1)^{\sigma(A,B)}e_{A\oplus B},$$  the twist function $\sigma(A, B)$ has the explicit form, given by
       \begin{equation}\label{0502-2056} \sigma(A,B)=  \left\{
 \begin{aligned} & 0, \quad  & &AB=0,
\\
&1, \quad  & &A=B\neq 0,
\\
       & b_l+\sum\limits_{i=l}^{n-1}\phi(a_i,b_i),\quad  & & A\neq B, \quad AB\neq 0,\quad \deg A> \deg B,\\
       & 1+a_l+\sum\limits_{i=l}^{n-1}\phi(a_i,b_i),\quad & &A\neq B, \quad AB\neq 0, \quad  \deg A<\deg B,\\
       & a_l+\sum\limits_{i=l}^{n-1}\phi(a_i,b_i),\quad  & &A\neq B,  \quad AB\neq 0,\quad \deg A=\deg B,
            \end{aligned}  \right.
     \end{equation}
where $A$ and $B$ are given by \eqref{eq:binary-490-23} and \eqref{eq:binary-491-23}, respectively, and 
$\phi$ is defined by  \eqref{eq:binary-490-99}  and
$l$ is defined in  \eqref{eq:phi-939}.
\end{theorem}

\begin{proof}  We note that the twist function $\sigma(A,B)$ satisfies certain properties. For example, for any $ A \in \{0, 1, \ldots, 2^n-1\}$, we have
\begin{eqnarray*}
\sigma(0, A)=\sigma(A, 0)=0
\end{eqnarray*}
and
\begin{eqnarray*}
   \sigma(A, A)
   =\begin{cases} 0, \qquad & A=0;
   \\ \\
  1, \qquad & A\neq 0.\end{cases}
\end{eqnarray*}

Hence, to derive an explicit formula for the twist function, it suffices to focus on the case where $1\leqslant A\neq B\leqslant 2^n-1$. By definition, we have $$\sigma(A, B)=\sigma(B, A)+1$$ in this case. Thus, we can further narrow down our analysis to the condition where $\deg A \geqslant \deg B$, subject to the constraints $1\leqslant A\neq B\leqslant 2^n-1$.  That is, we can  restrict  our attention to the scenario where  
\begin{eqnarray} \label{assump: twist-214} 1\leqslant A\neq B\leqslant 2^n-1, \qquad  \deg A \geqslant  \deg B. \end{eqnarray}

Henceforth, we will represent the integers using their binary expansions:
 \begin{eqnarray}\label{eq:binary-490} A&=&a_02^0+\dots+a_{n-1}2^{n-1},
 \\ \label{eq:binary-491}
     B&=&b_02^0+\dots+b_{n-1}2^{n-1}.
     \end{eqnarray}
All coefficients $a_j, b_j$  are either 0 or 1.

We define the degree of $A$ as the coefficient of the leading term in its binary representation. More precisely, we write $A$ in binary form as in \eqref{eq:binary-490} and set
   $$\deg\, A:=\min\{i\geqslant 0: a_i\neq 0   \}.$$

Assuming $1\leqslant A-a_{n-1}2^{n-1},B-b_{n-1}2^{n-1}\leqslant 2^{n-1}-1$, we have
 $$\sigma(A-a_{n-1}2^{n-1},B-b_{n-1}2^{n-1})+\sigma(B-b_{n-1}2^{n-1},A-a_{n-1}2^{n-1}=1.$$
 By \Cref{0504-1816}, we have
 \begin{equation}\label{0516-1503}
   \begin{aligned}
     \sigma(A,B)&=\sigma(A-a_{n-1}2^{n-1},B-b_{n-1}2^{n-1})(1+b_{n-1})+\sigma(B-b_{n-1}2^{n-1},A-a_{n-1}2^{n-1})b_{n-1}\\&\quad+\sigma(B-b_{n-1}2^{n-1},B-b_{n-1}2^{n-1})a_{n-1}+a_{n-1}b_{n-1}
     \\&=\left(\sigma(A-a_{n-1}2^{n-1},B-b_{n-1}2^{n-1})+\sigma(B-b_{n-1}2^{n-1},A-a_{n-1}2^{n-1})\right)b_{n-1}\\&\quad+\sigma(A-a_{n-1}2^{n-1},B-b_{n-1}2^{n-1})+a_{n-1}+a_{n-1}b_{n-1}
     \\&=\sigma(A-a_{n-1}2^{n-1},B-b_{n-1}2^{n-1})+b_{n-1}+a_{n-1}+a_{n-1}b_{n-1}
     \\&=\sigma(A-a_{n-1}2^{n-1},B-b_{n-1}2^{n-1})+\phi(a_{n-1},b_{n-1}).
   \end{aligned}
 \end{equation}

 By \eqref{0516-1503},
\begin{equation}\label{0516-1512}
  \begin{aligned}
    \sigma(A,B)&=\sigma(A-a_{n-1}2^{n-1},B-b_{n-1}2^{n-1})+\phi(a_{n-1},b_{n-1})
    \\&=\dots
    \\&=\sigma(A_l,B_l)+\sum\limits_{i=l+1}^{n-1}\phi(a_i,b_i).
  \end{aligned}
\end{equation}

Denote $$A_l=a_02^0+\dots+a_l2^l,$$
       $$B_l=b_02^0+\dots+b_l2^l.$$

 The proof is split into three cases.

 \bigskip
 {\bf Case 1}. $\deg A>\deg B.$
 \bigskip

In this particular case,   $l$ is equal  to    $\deg A$. As such, it follows that $A_l$ can be expressed as $2^l$.

Since
\begin{equation}
  e_{A_l}e_{B_l}=e_{2^l}e_{B_{l-1}+b_l2^l}=\left\{\begin{array}{ccc}
  -e_{B_{l-1}+2^l},&\mbox{if} & b_l=0,\\\\
  e_{B_{l-1}}, &\mbox{if}&  b_l=1,
 \end{array}\right.
\end{equation}
it follows that 
\begin{equation}\label{0516-1510}
  \sigma(A_l,B_l)=1+b_l.
\end{equation}

By substituting   \eqref{0516-1510} into  \eqref{0516-1512} and utilizing the fact that $1=\phi(a_l,b_l)$ as shown in \Cref{0515-1832}, it follows that 
$$\sigma(A,B)=1+b_l+\sum\limits_{i=l+1}^{n-1}\phi(a_i,b_i)=b_l+\sum\limits_{i=l}^{n-1}\phi(a_i,b_i).$$

\bigskip
{\bf Case 2}. $\deg A<\deg B.$
\bigskip

In this case,  we find that $l(A,B)$ and $l(B,A)$ are both equal to  $\deg B$. Furthermore, we know that $$\sigma(A,B)=1+\sigma(B,A).$$ As a result, we can apply the already established result from {\bf Case 1}. Therefore, we obtain the following expression:
\begin{equation*}
  \begin{aligned}
    \sigma(A,B)=1+\sigma(B,A)&=1+a_l+\sum\limits_{i=l(B,A)}^{n-1}\phi(b_i,a_i)
    \\&=1+a_l+\sum\limits_{i=l}^{n-1}\phi(a_i,b_i)
  \end{aligned}
\end{equation*}

\bigskip
{\bf Case 3}. $\deg A=\deg B.$

In this case, we have  
\begin{eqnarray*}l&=&\deg (A\oplus B), 
\\
A_{l-1}&=&B_{l-1}\neq 0, 
\\
a_l+b_l&=&1. 
\end{eqnarray*}

Since $$e_{A_l}e_{B_l}=e_{A_{l-1}+a_l2^l}e_{A_{l-1}+b_l2^l}=\left\{\begin{array}{ccc}
  e_{2^l},&\mbox{if}& a_l=1,b_l=0,\\\\
  -e_{2^l},&\mbox{if}& a_l=0,b_l=1,
\end{array}\right.$$
we have
\begin{equation}\label{0516-1539}
  \sigma(A_l,B_l)=1+a_l.
\end{equation}

We can substitute   \ eqref{0516-1539} into \eqref{0516-1512}. We observe that $1=\phi(a_l,b_l)$ and utilize the fact shown in \Cref{0515-1832}. This leads to the following efact: 
$$\sigma(A,B)=1+a_l+\sum\limits_{i={l+1}}^{n-1}\phi(a_i,b_i)=a_l+\sum\limits_{i={l}}^{n-1}\phi(a_i,b_i).$$
This finishes the proof.
\end{proof}

\begin{remark}
We observe from the proof of Theorem \ref{0515-1855} that it is sufficient to determine the value of $\sigma(A,B)$ only for the case where $1\leqslant A\neq B\leqslant 2^n-1$, as the function $\sigma(A,B)$ is trivially known for the remaining cases. Additionally, the relation between $\sigma(A,B)$ and $\sigma(B,A)$ is already established for this case.

Suppose that
$$1\leqslant A\neq B\leqslant 2^n-1, \quad
\mbox{and}\quad
\deg A\geqslant \deg B.$$
Then, we can express the twist function in a unified manner as follows:
\begin{equation}
\sigma(A,B)= \delta_{\deg A, \deg B}+b_l+\sum\limits_{i=l}^{n-1}\phi(a_i,b_i).
\end{equation}
\end{remark}

\begin{example} We take the example in the specific case that $n=3$.
When $n=3$, we have $2^n-1=7$, and the binary decompositions of $A$ and $B$ are as follows:
$$A=a_02^0+a_12^1+a_22^2,$$
$$B=b_02^0+b_12^1+b_22^2.$$
We assume   $A\neq B$ and $AB\neq 0$.  According to Theorem  \ref{0515-1855}, we can take the twist function as:
\begin{equation}\label{eq: twist-function-n=3}
\sigma(A,B) = \begin{cases}
b_2+\phi(a_2,b_2)+\phi(a_1,b_1)+\phi(a_0,b_0) & \text{if } \deg A > \deg B \\
a_2+\phi(a_2,b_2)+\phi(a_1,b_1)+\phi(a_0,b_0)+1 & \text{if } \deg A < \deg B \\
a_2+\phi(a_2,b_2)+\phi(a_1,b_1)+\phi(a_0,b_0) & \text{if } \deg A = \deg B
\end{cases}
\end{equation}

For example, suppose we take $A=5$ and $B=6$, then their binary decompositions are:
$$A=1\cdot 2^0+0\cdot 2^1+1\cdot 2^2,$$
$$B=0\cdot 2^0+1\cdot 2^1+1\cdot 2^2.$$
Therefore, we have $\deg A=\deg B=2$, and $A\oplus B=3$ with $\deg(A\oplus B)=1$. Thus, we have $l=2$, and we can use the formula in \eqref{eq: twist-function-n=3} to compute the twist function:
$$\sigma(A,B) = a_2+\phi(a_2,b_2)+\phi(a_1,b_1)+\phi(a_0,b_0) = 1+1+1+2 = 5.$$
This result matches the value of $\sigma(A,B)$ that can be obtained directly from the multiplication table of the octonions. Therefore, Theorem   \ref{0515-1855} is correct in the case of $n=3$.
\end{example}

\section{Split Cayley-Dickson algebras}

In the following section, we focus on split Cayley-Dickson algebras and provide explicit formulas for their twist functions. We begin by defining an involution and regular involution on an algebra over a field $\mathbb F$ using the doubling product definition from \cite{S1954}. A Cayley-Dickson algebra is an algebra with a regular involution satisfying specific conditions. Furthermore, a split Cayley-Dickson algebra is a subtype of Cayley-Dickson algebra that has a split norm form. The formulas we present in this section are derived using the same approach as before.

\begin{definition}
Let $(\mathcal A,+,\cdot)$ be an algebra over a field $\mathbb F$. An involution $a\mapsto \bar a$ on $\mathcal A$ is an endomorphism of the vector space $\mathcal A$ such that for $a,b\in \mathcal A,$ we have $$\bar{\bar a}=a$$ and $$\overline{ab}=\bar b\bar a.$$

Let $1_{\mathcal A}$ be the identity of $\mathcal A.$ An involution $a\mapsto \bar a$ is called regular if for any $a\in \mathcal A,$ we have $$a+\bar a=t(a)1_{\mathcal A},$$ $$a\bar a=\bar a a=n(a)1_{\mathcal A},$$ where $t(a)\in \mathbb F$ is called the trace of $a$ and $n(a)\in \mathcal F$ is called the norm of $a$.
\end{definition}

\begin{definition}
 Let us assume that $\mathcal A$ is an algebra over a field $\mathbb F$ with a regular involution $a\mapsto \bar a$.
We define the algebra $\mathcal A\{\gamma\}$ to be the result of applying the Cayley-Dickson process to $\mathcal A$ using a parameter $\gamma\in \mathbb F$, where $\gamma \neq 0$. The elements of
 $\mathcal A\{\gamma\}$ are ordered pairs of elements of $\mathcal A$ and the operations are defined as follows:

$\bullet$ $\alpha (a,b)=(\alpha a,\alpha b)$;

$\bullet$  $(a,b)+(c,d)=(a+c,b+d)$;

$\bullet$ 
$(a,b)(c,d)=(ac+\gamma \bar d b,da+b\bar c)$;

 and
$$\overline{(a,b)}=(\bar a,-b).$$
Here, $a,b,c,d\in \mathcal A$ and $\alpha \in \mathbb F$.
\end{definition}

From this definition, we can express the Cayley-Dickson algebra $\mathbb D$ as:
$$\mathbb D=\mathbb D^{(n)}=\R\{\gamma_0,\dots,\gamma_{n-2},\gamma_{n-1}\}$$
where $\gamma_0=\dots=\gamma_{n-1}=-1$.

We can use induction to obtain the following expressions:
\begin{align}
\mathbb D^{(n+1)} &= \mathbb D (\gamma_n) = \R\{\gamma_0,\dots,\gamma_{n-2},\gamma_{n-1}, \gamma_{n} \}, \\
\mathbb D^{(n+1)}_s &= \mathbb D\{\gamma_n^s\} = \R\{\gamma_0,\dots,\gamma_{n-2},\gamma_{n-1},\gamma_n^s\},
\end{align}
where $\gamma_0=\dots=\gamma_{n-1}=\gamma_{n}=-1$ and $\gamma_{n}^s=1$.

Now, we define the twist function $$\sigma_s:\mathbb Z \times \mathbb Z \longrightarrow \mathbb Z_2$$ for the split Cayley-Dickson algebra $\mathbb D^{(n+1)}_s$, as follows:
 $$e_pe_q=(-1)^{\sigma_s(p,q)}e_{p\oplus q}$$
where $p,q=0,\cdots,2^{n+1}-1.$ This definition is from \cite{AM1999}.

We can prove the following recursion formula for the split function $\sigma_s$.

\begin{proposition}\label{0504-1810} For   any integers $0\leqslant A,B\leqslant 2^n-1$ and $a_n, b_n\in\{0, 1\}$, we have 
\begin{equation}\label{0504-1805}
\sigma_s(A+a_{n}2^{n},B+b_{n}2^{n})=\sigma_s(A,B)(1+b_{n})+\sigma_s(B,A)b_{n}+\sigma_s(B,B)a_{n}.
\end{equation}
\end{proposition}

\begin{proof}
We use the following facts:
\begin{align}
    e_A(e_Be_{2^n}) &= (e_Be_A)e_{2^n} \label{1825-1}\\
    (e_Ae_{2^n})e_B &= (e_A\bar e_B)e_{2^n} \label{1825-2}\\
    (e_Ae_{2^n})(e_Be_{2^n}) &= \bar e_B e_A \label{1825-3}
\end{align}

Fact \eqref{1825-1} implies that $$\sigma_s(A,B+2^n) = \sigma_s(B,A).$$ fact \eqref{1825-2} implies that $$\sigma_s(A+2^n,B) = \sigma_s(A,B)+\sigma_s(B,B).$$ and fact \eqref{1825-3} implies that $$\sigma_s(A+2^n,B+2^n) = \sigma_s(B,A)+\sigma_s(B,B).$$

Combining these facts, we get:
\begin{align*}
    \sigma_s(A+a_{n}2^{n},B+b_{n}2^{n}) &= \sigma_s(A+a_n2^n,B)+\sigma_s(B,b_n2^n)\\
    &= [\sigma_s(A,B)+\sigma_s(B,B)a_n]+\sigma_s(B,A)(1+b_n)\\
    &= \sigma_s(A,B)(1+b_n)+\sigma_s(B,A)b_n+\sigma_s(B,B)a_n,
\end{align*}
which proves the proposition.
\end{proof}

The explicit form of the twist function for the split Cayley-Dickson algebra can be reduced to that of the Cayley-Dickson algebra.

\begin{theorem}
  Let $0\leqslant A',B'\leqslant 2^{n+1}-1$ with binary decompositions
   \begin{eqnarray*}
   A^\prime=a_02^0+\dots+a_{n-1}2^{n-1}+a_n2^n,
   \\
     B^\prime=b_02^0+\dots+b_{n-1}2^{n-1}+b_n2^n.
     \end{eqnarray*}
   Suppose $\sigma_s$ is the twist function of the split Cayley-Dickson algebra $\mathbb D^{(n+1)}_s$ and $\sigma$ is the twist function of the Cayley-Dickson algebra $\mathbb D^{(n+1)}$. Then, we have
           \begin{equation}\label{0504-1724}
         \sigma_s(A^\prime,B^\prime)=\sigma(A^\prime,B^\prime)+a_{n}b_{n}.
     \end{equation}
\end{theorem}

\begin{proof}
Let $A^\prime = A+a_n2^n$ and $B^\prime = B+b_n2^n$. Using \Cref{0504-1810}, we have
\begin{align*}
\sigma_s(A^\prime,B^\prime) &= \sigma_s(A+a_{n}2^{n},B+b_{n}2^{n}) \\
&= \sigma_s(A,B)(1+b_{n})+\sigma_s(B,A)b_{n}+\sigma_s(B,B)a_{n} \\
&= \sigma(A,B)(1+b_{n})+\sigma(B,A)b_{n}+\sigma(B,B)a_{n} \\
&= a_nb_n+\sigma(A+2^n,B+2^n) \\
&= a_nb_n+\sigma(A^\prime,B^\prime),
\end{align*}
where the penultimate equality follows from \Cref{0504-1816}. Therefore, we have proved \eqref{0504-1724}.
\end{proof}

\bigskip\bigskip 

\end{document}